\documentclass[12pt]{article}

\usepackage{amsmath}
\usepackage{amsfonts}
\usepackage{amsthm}
\usepackage{graphics}
\newtheorem{theorem}{Theorem}[section]
\newtheorem{lemma}{Lemma}
\numberwithin{equation}{section}

\def\ds{\displaystyle}
\def \reff#1{(\ref{#1})}
\begin{document}
\title{Orthogonal polynomials on a bi-lattice \thanks{Research supported by FWO research grant G.0427.09, K.U.Leuven research grant OT/08/033
and the Belgian Interuniversity Attraction Poles Programme P6/02.}}
\author{Christophe Smet and Walter Van Assche \\ Katholieke Universiteit Leuven, Belgium}
\maketitle
\begin{abstract}
We investigate generalizations of the Charlier and the Meixner polynomials on the lattice $\mathbb{N}$ and on the shifted
lattice $\mathbb{N} + 1-\beta$. We combine both lattices to obtain the bi-lattice $\mathbb{N} \cup (\mathbb{N} + 1-\beta)$
and show that the orthogonal polynomials on this bi-lattice have recurrence coefficients which satisfy a non-linear system
of recurrence equations, which we can identify as a limiting case of an (asymmetric) discrete Painlev\'e equation.
\end{abstract}

\begin{description}
  \item[2010 Mathematics Subject Classification] 33C47, 42C05, 34M55, 65Q30. 
  \item[Keywords] Discrete orthogonal polynomials, recurrence coefficients, discrete Painlev\'e equations. 
\end{description}      

\section{Introduction}
The classical orthogonal polynomials (Jacobi, Laguerre and Hermite polynomials) are orthogonal polynomials
on (an interval) of the real line with a weight function that satisfies a first order differential equation
(the so-called Pearson equation) of the form
\[     [\sigma(x) w(x)]' = \tau(x) w(x), \]
where $\sigma$ is a polynomial of degree at most 2 and $\tau$ is a polynomial of degree 1.
This Pearson equation allows to find many useful properties of these polynomials.
There are other families of orthogonal polynomials in the Askey table \cite{KLS} which live on a linear lattice
rather than on an interval. These are the Charlier polynomials (\cite[\S 9.14 on p.~247]{KLS}, \cite[Chapter VI, \S 1]{Chihara})
\[   \sum_{k=0}^\infty C_n(k;a) C_m(k;a) \frac{a^k}{k!} = a^{-n} e^a n! \delta_{n,m}, \qquad a > 0, \]
which are orthogonal on the lattice $\mathbb{N}$ with respect to the Poisson distribution, and the Meixner polynomials
(\cite[\S 9.10 on p.~234]{KLS}, \cite[Chapter VI, \S 3]{Chihara} who calls them Meixner polynomials of the first kind)
\[    \sum_{k=0}^\infty M_n(k;\beta,c) M_m(k;\beta,c) \frac{(\beta)_k c^k}{k!} = \frac{c^{-n}n!}{(\beta)_n(1-c)^\beta}  \delta_{m,n},
\qquad \beta > 0,\ 0 < c < 1, \]
which are orthogonal on $\mathbb{N}$ with respect to the negative binomial (or Pascal) distribution. The special case $\beta=-N$
and $c = p/(p-1)$, with $N \in \mathbb{N}$ and $0 < p <1$ are the Krawtchouk polynomials which are orthogonal on the integers
$\{0,1,2,\ldots, N\}$ with respect to the binomial distribution. The Hahn polynomials are another family of polynomials
which are orthogonal on the integers $\{0,1,2,\ldots, N\}$.
Instead of the differential operator it is better for these lattice polynomials to work with the difference operators
\[     \nabla f(x) = f(x)-f(x-1), \qquad \textrm{(backward difference)}  \]
and
\[   \Delta f(x) = f(x+1)-f(x), \qquad \textrm{(forward difference)}. \]
The weight $w_k = a^k/k!$ for Charlier polynomials now satisfies the Pearson equation
\[   \nabla w(x) = \left( 1 - \frac{x}{a} \right) w(x)  \]
if we define the weight function $w(x) = a^x/\Gamma(x+1)$, so that $w_k = w(k)$. This is of the form
$\nabla[\sigma(x)w(x)] = \tau(x) w(x)$ with $\sigma=1$ and $\tau$ a polynomial of degree 1. The weight for Meixner
polynomials $w_k= (\beta)_k c^k/k!$ can be written as $w_k=w(k)$ using the weight function $w(x) = \Gamma(\beta+x)c^x/(\Gamma(\beta)\Gamma(x+1))$
and this satisfies the Pearson equation
\[  \nabla [(\beta+x)w(x)] = \left( \beta+x-\frac{x}{c} \right) w(x), \]
so that $\sigma(x)=\beta+x$ and $\tau$ is a polynomial of degree 1. This Pearson equation (with a difference operator) allows to find many
properties of these orthogonal polynomials on the lattice $\mathbb{N}$. In particular one can find the three term recurrence relation
for these polynomials explicitly:
\[   -x C_n(x;a) = a C_{n+1}(x;a) - (n+a) C_n(x;a) + n C_{n-1}(x;a)  \]
and
\[  (c-1)xM_n(x;\beta,c) = c(n+\beta)M_{n+1}(x;\beta,c) - [n+(n+\beta)c]M_n(x;\beta,c) + nM_{n-1}(x;\beta,c). \]
Often it is more convenient to work with the monic polynomials $P_n(x) = x^n + \cdots$, for which the recurrence relation is
\[    xP_n(x) = P_{n+1}(x) + b_n P_n(x) + a_n^2 P_{n-1}(x), \]
or with the orthonormal polynomials $p_n(x) = \gamma_n x^n + \cdots$, with
\[   \frac{1}{\gamma_n^2} = \sum_{k=0}^\infty P_n^2(k) w_k, \]
for which the recurrence relation is
\[    xp_n(x) = a_{n+1} p_{n+1}(x) + b_n p_n(x) + a_n p_{n-1}(x). \]
The recurrence coefficients for Charlier polynomials then can be seen to be $a_n^2 = na$ and $b_n = n+a$ for $n \in \mathbb{N}$, and for
Meixner polynomials one has
\[   a_n^2 = \frac{n(n+\beta-1)c}{(1-c)^2}, \quad b_n = \frac{n+(n+\beta)c}{1-c},\qquad n \in \mathbb{N}. \]

In this paper we will investigate generalizations of the Charlier and the Meixner polynomials. In Section 2.1 we will use the weight
function $w(x) = \Gamma(\beta) a^x/(\Gamma(\beta+x) \Gamma(x+1))$ which gives the weights $w_k=w(k)= a^k/((\beta)_k k!)$ on the lattice
$\mathbb{N}$. In Theorem \ref{thm:1} we will find a system of non-linear recurrence equations for the recurrence coefficients $a_n^2$ and $b_n$
and we can identify this system as a limiting case of a discrete Painlev\'e IV equation. In Section 2.2 we will use this weight
on the shifted lattice $\mathbb{N} + 1-\beta$ and in Theorem \ref{thm:2} we will show that the recurrence coefficients satisfy the
same non-linear system of recurrence equations, but with a different initial condition for $b_0$. In Section 2.3 we combine both lattices
to obtain the bi-lattice $\mathbb{N} \cup (\mathbb{N} + 1-\beta)$. The orthogonality measure is now a linear combination of the
measures on $\mathbb{N}$ and $\mathbb{N}+ 1-\beta$ and surprisingly, Theorem \ref{thm:3} shows that the recurrence coefficients again
satisfy the same non-linear system of recurrence equations, but with a different initial condition for $b_0$. In all these cases, the initial conditions
for $b_0$ are given by ratios of modified Bessel functions.
In Section 3 we deal with a similar generalization of the Meixner weight and we investigate the weight function
$w(x) = \Gamma(\beta) \Gamma(\gamma+x)a^x/(\Gamma(\gamma)\Gamma(\beta+x) \Gamma(x+1))$. In Section 3.1 we use this weight on the lattice $\mathbb{N}$
where it gives the weights $w_k=w(k) = (\gamma)_k a^k/((\beta)_k k!)$. In Theorem \ref{GenMeixThmOriginalLattice} we will show that
the recurrence coefficients of the corresponding orthogonal polynomials satisfy another non-linear system of recurrence equations, which we can
identify as a limiting case of an asymmetric discrete Painlev\'e IV equation. In Section 3.2 we use the same weight function but now on the
shifted lattice $\mathbb{N} + 1-\beta$. The recurrence coefficients of the corresponding orthogonal polynomials satisfy the same non-linear
system of recurrence equations, but with a different initial value for $b_0$. Finally, in Section 3.3 we study the orthogonal polynomials
on the bi-lattice $\mathbb{N} \cup (\mathbb{N}+1-\beta)$ with respect to a linear combination of the measures on $\mathbb{N}$ and $\mathbb{N} + 1-\beta$,
and again the recurrence coefficients satisfy the same system of non-linear recurrence equations, but with a different initial value for $b_0$.
The initial conditions in Section 3 are all given by ratios of confluent hypergeometric functions.

We hope this paper is of interest to people studying orthogonal polynomials because we introduce some new families of orthogonal polynomials
which seem to have nice properties. The asymptotic behavior of the recurrence coefficients may be of interest since in many cases one can observe an
oscillating behavior, reminiscent of orthogonal polynomials on two intervals. The paper should also be of interest to people studying discrete
Painlev\'e equations since we are giving two systems of discrete Painlev\'e equations and a one-parameter family of initial conditions
for which the solution turns out to be in terms of recurrence coefficients of orthogonal polynomials on the real line. This means that these particular
initial values do not lead to singularities and in fact have nice positivity properties, such as $a_n >0$. The fact that the polynomials are
orthogonal on a (bi-)lattice also gives some properties for partial sums of the $b_n$ (sums of the zeros).
The cases where the bi-lattice reduces to a single lattice may well be corresponding to unique solutions of the discrete Painlev\'e equations with a
prescribed asymptotic behavior for the $b_n$ as $n \to \infty$.

\section{A generalized Charlier weight}
\subsection{The lattice $\mathbb{N}$}
We consider the discrete weights
\[    w_k = \frac{a^k}{(\beta)_k k!}, \qquad k \in \mathbb{N} = \{0,1,2,3,\ldots \}
\] with $a>0$.
We can write these in terms of the weight function
\begin{equation} \label{eq:w}
    w(x) = \frac{\Gamma(\beta) a^x}{\Gamma(\beta+x) \Gamma(x+1)},
\end{equation}
which is now a function of $x \in \mathbb{C}$ vanishing whenever $x$ is a pole of
$\Gamma(x+1)$ or a pole of $\Gamma(\beta+x)$, i.e., for $x = -1,-2,-3,\ldots$ and $x=-\beta, -\beta-1,
-\beta-2,\ldots$. This weight satisfies the Pearson equation
\begin{equation}  \label{eq:Pearson}
     \nabla w(x) := w(x)-w(x-1) = \frac{ a-x(\beta-1)-x^2}{a} w(x).
\end{equation}
With this weight we can introduce the inner product
\[  \langle f,g \rangle = \sum_{k=0}^\infty  f(k)g(k) w_k   \]
which has positive weights for every $\beta > 0$. We denote the corresponding monic
orthogonal polynomials by $P_n(x;a,\beta)$:
\begin{equation}   \label{eq:Pbeta}
   \sum_{k=0}^\infty P_n(k;a,\beta)P_m(k;a,\beta) \frac{a^k}{(\beta)_k k!} = 0, \qquad n \neq m.
\end{equation}
Our main interest is to find the recurrence coefficients in the three term recurrence
relation
\begin{equation}  \label{eq:Prec}
     xP_n(x;a,\beta) = P_{n+1}(x;a,\beta) + b_n P_n(x;a,\beta) + a_n^2 P_{n-1}(x;a,\beta),
\end{equation}
with initial conditions $P_{-1}=0$ and $P_0=1$. We will obtain
recursive relations for these recurrence coefficients. We sometimes
just write $P_n(x)$ for the polynomials $P_n(x;a,\beta)$ to simplify
the notation. Observe that for $a = \beta c$ and $\beta \to \infty$
we get the weights $w_k = c^k/k!$ so that
\[   \lim_{\beta \to \infty} P_n(x;\beta c,\beta) = \hat{C}_n(x;c)  \]
are the monic Charlier polynomials.

The Pearson equation and summation by parts gives the following structure relation for the
orthogonal polynomials
\begin{lemma}
The monic orthogonal polynomials given by (\ref{eq:Pbeta}) satisfy the relation
\begin{equation}  \label{eq:Pstruc}
      \Delta P_n(x) := P_n(x+1)-P_n(x) = nP_{n-1}(x) + B_nP_{n-2}(x),
\end{equation}
for some sequence $B_n$ of real numbers.
\end{lemma}

\begin{proof}
We can always write the polynomial $P_n(x+1)-P_n(x)$ (which is of degree $n-1$) in a Fourier series using the
orthogonal polynomials and hence
\[    P_n(x+1)-P_n(x) = \sum_{k=0}^{n-1} A_{k,n} P_k(x),  \]
where the Fourier coefficients are given by
\[      A_{k,n} = \frac{\langle \Delta P_n, P_k \rangle}{\langle P_k, P_k \rangle} .  \]
Recall summation by parts: if $b_{-1}=0$ then
\[    \sum_{k=0}^\infty (\Delta a_k) b_k = - \sum_{k=0}^\infty a_k \nabla b_k   \]
whenever $a$ and $b$ are in $\ell_2$. If we apply this, then
\begin{eqnarray*}   \langle P_k,P_k \rangle A_{k,n} &=& \sum_{j=0}^\infty (\Delta P_n(j)) P_k(j) w_j \\
                     &=& - \sum_{j=0}^\infty P_n(j) \nabla (P_k(j)w_j) \\
                     &=& - \sum_{j=0}^\infty P_n(j) w_j \nabla P_k(j)
                         - \sum_{j=0}^\infty P_n(j) P_k(j-1) \nabla
                         w_j.
\end{eqnarray*}
The first sum on the right is zero by orthogonality since $\nabla
P_k$ is a polynomial of degree $k-1 < n$. For the second sum we use
the Pearson equation (\ref{eq:Pearson}) to find
\[    \langle P_k,P_k \rangle A_{k,n} = - \frac{1}{a} \sum_{j=0}^\infty P_n(j)P_{k}(j-1) (a-j(\beta-1)-j^2) w_j. \]
This sum is zero by orthogonality whenever $k+2 < n$. Hence only
$A_{n-1,n}$ and $A_{n-2,n}$ can be non-zero. By comparing the
leading coefficients, we see that $A_{n-1,n} = n$. If we call
$A_{n-2,n}=B_n$, then the required formula follows.
\end{proof}

Some of these coefficients are useful in later computations. If we define
\[   \frac{1}{\gamma_n^2} = \langle P_n, P_n \rangle, \]
then it is not difficult to find from the recurrence relation (\ref{eq:Prec}) (taking the inner product with
$P_{n-1}$) the well-known relation
\[   a_n^2 = \frac{\gamma_{n-1}^2}{\gamma_n^2}. \]
Furthermore, with the notation in the proof of the lemma
\[   1 = A_{0,1} =  -\frac{\gamma_0^2}{a} \sum_{j=0}^\infty P_1(j) (a-j(\beta-1)-j^2) w_j. \]
By orthogonality we have
\[   \sum_{j=0}^\infty P_1(j)w_j = 0.  \]
If we use the recurrence relation then
\[   \sum_{j=0}^\infty j P_1(j) w_j = \sum_{j=0}^\infty [P_2(j)+b_1P_1(j)+a_1^2P_0(j)]w_j = a_1^2/\gamma_0^2. \]
If we use the recurrence relation twice, then
\begin{eqnarray*}
       \sum_{j=0}^\infty j^2 P_1(j) w_j &=& \sum_{j=0}^\infty [b_1(P_2(j) + b_1P_1(j)+a_1^2P_0(j))+a_1^2(P_1(j)+
b_0P_0(j))]w_j \\
    &=& a_1^2(b_1+b_0)/\gamma_0^2.
\end{eqnarray*}
Hence, combining these results, we have
\begin{equation}  \label{eq:A1}
     1 = \frac{a_1^2}{a} (b_1 + b_0 + \beta-1).
\end{equation}
In a similar way we can compute $B_2$:
\[    B_2 = A_{0,2} = -\frac{\gamma_0^2}{a} \sum_{j=0}^\infty P_2(j) (a-j(\beta-1)-j^2) w_j. \]
By orthogonality
\[   \sum_{j=0}^\infty P_2(j)w_j = 0 = \sum_{j=0}^\infty j P_2(j) w_j \]
and using the recurrence relation we get
\begin{eqnarray*}
    \sum_{j=0}^\infty j^2 P_2(j)w_j &=& \sum_{j=0}^\infty j(P_3(j)+b_2P_2(j)+a_2^2 P_1(j)) w_j \\
                         &=& a_2^2a_1^2/\gamma_0^2,
\end{eqnarray*}
so that
\begin{equation}   \label{eq:B2}
     B_2 = \frac{a_1^2a_2^2}{a}.
\end{equation}

We now are able to derive the recursive relations for the recurrence coefficients.

\begin{theorem}  \label{thm:1}
The recurrence coefficients for the orthogonal polynomials  defined by
(\ref{eq:Pbeta}) satisfy
\begin{eqnarray}
   b_n+b_{n-1}-n+\beta &=& \frac{an}{a_n^2},  \label{eq:rec1}  \\
  (a_{n+1}^2-a)(a_n^2-a) &=& a (b_n-n) (b_n-n+\beta-1), \label{eq:rec3}
\end{eqnarray}
with initial conditions
\[   a_0^2=0, \qquad   b_0 = \frac{\sqrt{a} I_\beta(2\sqrt{a})}{I_{\beta-1}(2\sqrt{a})}, \]
where $I_\nu$ is the modified Bessel function
\[   I_\nu(z) = \sum_{k=0}^\infty \frac{(z/2)^{2k+\nu}}{k! \Gamma(k+\nu+1)}.   \]
\end{theorem}

\begin{proof}
On one hand we have the three-term recurrence relation (\ref{eq:Prec})
\[   kP_n(k) = P_{n+1}(k) + b_n P_n(k) + a_n^2 P_{n-1}(k), \]
and on the other hand we have the structure relation (\ref{eq:Pstruc}). The compatibility between these two equations
will result in the desired equations. Take the forward difference of the three-term recurrence relation
to find
\[  (k+1) \Delta (P_n(k)) + P_n(k) = \Delta P_{n+1}(k) + b_n \Delta P_{n}(k) + a_n^2 \Delta P_{n-1}(k). \]
Then use the structure relation (\ref{eq:Pstruc}) to find
\begin{multline*}   (k+1) \left( nP_{n-1}(k) + B_nP_{n-2}(k) \right) + P_n(k) \\
  =    (n+1) P_n(k) + B_{n+1} P_{n-1}(k) + b_n \left( nP_{n-1}(k) + B_n P_{n-2}(k)
  \right)\\
   + a_n^2 \left( (n-1)P_{n-2}(k) + B_{n-1} P_{n-3}(k) \right).
\end{multline*}
Finally, use the three-term recurrence relation to work out $kP_{n-1}(k)$ and $kP_{n-2}(k)$ on the left hand side.
This gives a linear identity involving $P_n(k), P_{n-1}(k), P_{n-2}(k), P_{n-3}(k)$, and since these four polynomials (in the variable $k$)
are linearly independent, the coefficients should all be zero. This gives the four equations
\begin{eqnarray}
  n+1 &=& n+1,   \label{eq:comp1}\\
  n+nb_{n-1} + B_n &=& B_{n+1} + nb_n, \label{eq:comp2} \\
  B_n + na_{n-1}^2 + B_n b_{n-2} &=& B_n b_n + (n-1) a_n^2,  \label{eq:comp3} \\
    B_n a_{n-2}^2 &=& a_n^2 B_{n-1}.   \label{eq:comp4}
\end{eqnarray}
Clearly (\ref{eq:comp1}) is always satisfied and (\ref{eq:comp4}) readily gives
\[ \frac{B_n}{a_n^2a_{n-1}^2} = \frac{B_{n-1}}{a_{n-1}^2a_{n-2}^2}, \]
so that $B_n =a_n^2a_{n-1}^2 B_2/(a_2^2a_1^2)$, and hence (\ref{eq:B2}) gives
\begin{equation}  \label{eq:Bn}
    B_n = \frac{a_n^2a_{n-1}^2}{a}, \qquad n \geq 2.
\end{equation}
Use this in (\ref{eq:comp2}) to find
\begin{equation}  \label{eq:rec2}
    na \left( b_n-b_{n-1}-1 \right) = a_n^2 (a_{n-1}^2-a_{n+1}^2),
\end{equation}
and (\ref{eq:comp3}) becomes
\begin{equation*}
    \frac{1}{a} ( b_n - b_{n-2} -1 ) = \frac{n}{a_n^2} - \frac{n-1}{a_{n-1}^2}.
\end{equation*}
Summing the latter starting from $n=2$ gives
\begin{equation*}
  \frac{1}{a} (b_n+b_{n-1}-n+1) - \frac{1}{a}(b_1+b_0) = \frac{n}{a_n^2} - \frac{1}{a_1^2}.
\end{equation*}
Now use (\ref{eq:A1}) to find (\ref{eq:rec1}).

If we use \eqref{eq:rec1} in \eqref{eq:rec2} and put $b_n = n + d_n$, then we find
\[ (d_k+d_{k-1}+k+\beta-1)(d_{k-1}-d_{k}) = a_{k+1}^2-a_{k-1}^2.   \]
Summing from $k=1$ to $n$ gives (use the telescoping property and summation by parts)
\begin{equation}  \label{eq:elim1}
   -d_n^2 +d_0^2 + \sum_{k=0}^{n-1} d_k - nd_n - (\beta-1) (d_n-d_0) = a_{n+1}^2+a_n^{2}-a_1^2,
 \end{equation}
where we used the initial condition $a_0^2=0$. On the other hand, \eqref{eq:rec2} is equivalent with
\[   ak(d_{k-1}-d_k) = a_k^2a_{k+1}^2 - a_{k-1}^2a_k^2. \]
Summing from $k=1$ to $n$ now gives
\begin{equation}   \label{eq:elim2}
   a \sum_{k=0}^{n-1} d_k - and_n = a_n^2a_{n+1}^2.
\end{equation}If we use \eqref{eq:elim2} in \eqref{eq:elim1} then we find
\begin{equation}  \label{eq:rec12}
      -d_n^2+d_0^2 + \frac{a_n^2a_{n+1}^2}{a} -(\beta-1)(d_n-d_0) = a_{n+1}^2+a_n^2-a_1^2.
\end{equation}
The initial values $d_0=b_0$ and $a_1^2$ are given by
\[    b_0 = \frac{m_1}{m_0}, \quad a_1^2 = \frac{m_2}{m_0} - \left( \frac{m_1}{m_0} \right)^2, \]
where $m_j$ are the moments
\[   m_j = \sum_{k=0}^\infty k^j w_k.  \]
A simple calculation gives
\[   m_0 = \frac{\Gamma(\beta)}{(\sqrt{a})^{\beta-1}} I_{\beta-1}(2\sqrt{a}), \quad
     m_1 = \frac{\Gamma(\beta)}{(\sqrt{a})^{\beta-2}} I_{\beta}(2\sqrt{a}), \]
where $I_\beta$ and $I_{\beta-1}$ are modified Bessel functions. This gives
\begin{equation}   \label{eq:b0}
    b_0 = \sqrt{a} \frac{I_\beta(2\sqrt{a})}{I_{\beta-1}(2\sqrt{a})}.
\end{equation}The Pearson equation \eqref{eq:Pearson} gives
\[  m_2 = \sum_{k=0}^\infty k^2 w_k = \sum_{k=0}^\infty [a-k(\beta-1)] w_k - a \sum_{k=0}^\infty [w_k-w_{k-1}], \]
which readily gives $m_2 = a m_0 - (\beta-1) m_1$. Observe that this gives
\[       d_0^2+(\beta-1)d_0 + a_1^2 = a  \]
which simplifies the recurrence relation \eqref{eq:rec12} to
\[      a_{n+1}^2+a_n^2 + d_n^2 - \frac{a_n^2a_{n+1}^2}{a} + (\beta-1)d_n - a = 0,   \]
which is equivalent with (\ref{eq:rec3}).
\end{proof}

The case $\beta = 1$ was already considered in \cite{mama} (see also \cite[\S 4.2]{WVA}). Equation (\ref{eq:rec3}) then simplifies to
\begin{equation}  \label{eq:rec31}
    (a_{n+1}^2-a)(a_n^2-a) = (b_n-n)^2,
\end{equation}
so that $a_n^2-a$ and $a_{n+1}^2-a$ have the same sign. The boundary condition $a_0^2=0$ thus implies that $a_n^2-a < 0$ for all $n$, and we can write
$a_n^2 -a = -ac_n^2$, for some new sequence $(c_n)_{n\in \mathbb{N}}$. Then $a_n^2 = a(1-c_n^2)$ so that $c_n^2 < 1$ for $n \geq 1$. This
still leaves two choices for the sign of $c_n$. Taking square roots in (\ref{eq:rec31}) gives $b_n = n + a c_n c_{n+1}$, where we choose
$c_0=1$ and we recursively take the sign of $c_{n+1}$ equal to the sign of $(b_n-n)/c_n$. Inserting these formulas for $a_n^2$ and $b_n$ into
(\ref{eq:rec1}) with $\beta = 1$ gives $a(c_{n+1}+c_{n-1}) = n/(1-c_n^2)$, which is the discrete Painlev\'e II equation (\cite{GR},
\cite[Appendix A.1]{WVA}).

The situation is different for $\beta \neq 1$, since (\ref{eq:rec3}) now is a quadratic equation in $b_n$ and it is not a priori
clear which of the two roots one should choose. Equations (\ref{eq:rec1})--(\ref{eq:rec3}) are a limiting case of a discrete Painlev\'e IV
equation: take the first dP$_{\textrm{IV}}$ in \cite[p.~723]{WVA}, or equivalently the second discrete Painlev\'e equation for $D_4^c$ in
\cite[p.~297]{GR}
\begin{eqnarray*}
    x_{n+1}x_n &=& \frac{(y_n-z_n)^2-A}{y_n^2-B} \\
      y_n+y_{n-1} &=& \frac{\zeta_n-C}{1+Dx_n} + \frac{\zeta_n+C}{1+x_n/D}
\end{eqnarray*}
where $z_n=z_0+n\delta$ and $\zeta_n = z_n-\delta/2$. If we put $x_n = iX_n/\sqrt{aB}$ and let $B \to \infty$ then
this gives for the first equation
\[    X_{n+1}X_n = a\left( (y_n-z_n)^2-A \right), \]
and if we take $iD = \sqrt{B/a}$ and let $B \to \infty$ then we also get
\[    y_n+y_{n-1} = \frac{\zeta_n-C}{1+X_n/a} + \zeta_n+C. \]
The parameters $A=(\beta-1)^2/4$, $C= -\beta/2$, $z_n = n-(\beta-1)/2$ then give the discrete equations
(\ref{eq:rec1})--(\ref{eq:rec3}) for $X_n=a_n^2-a$ and $y_n=b_n$.

\subsection{The shifted lattice $\mathbb{N}+1-\beta$}
We can consider the weight $w$ in (\ref{eq:w}) also on the shifted lattice $\mathbb{N}+1-\beta = \{1-\beta, 2-\beta, 3-\beta, \ldots\}$,
where
\[    v_k := w(k+1-\beta) =  \frac{\Gamma(\beta)a^{1-\beta}}{\Gamma(2-\beta)} \frac{a^k}{k!(2-\beta)_k}  , \qquad k\in\mathbb{N}.  \]
The weights $(v_k)_{k\in \mathbb{N}}$ are therefore of the same form as in the previous section (up to a real factor) but with $\beta$
replaced by $2-\beta$.
The corresponding monic orthogonal polynomials $Q_n$ on the shifted lattice $\mathbb{N}+1-\beta$ satisfy
\begin{equation}  \label{eq:Qbeta}
    \sum_{k=0}^\infty Q_n(k+1-\beta)Q_m(k+1-\beta) v_k = 0, \qquad n \neq m,
\end{equation}
and they are simply the polynomials $P_n(\cdot;a,2-\beta)$ but shifted over a distance $1-\beta$:
\begin{equation}  \label{eq:PQ}
    Q_n(x) = P_n(x+\beta-1;a,2-\beta).
\end{equation}
These are orthogonal polynomials with a positive measure whenever $\beta < 2$. The remarkable fact is that the recurrence coefficients
in the three-term recurrence relation
\begin{equation} \label{eq:Qrec}
     xQ_n(x) = Q_{n+1}(x) + \hat{b}_n Q_n(x) + \hat{a}_n^2 Q_{n-1}(x)
\end{equation}
satisfy the same system of non-linear recurrence relations (discrete Painlev\'e equations) as in Theorem \ref{thm:1}, but with
a different initial condition.

\begin{theorem}  \label{thm:2}
The recurrence coefficients for the orthogonal polynomials defined by (\ref{eq:Qbeta}) satisfy the system of equations
\begin{eqnarray*}
   \hat{b}_n+\hat{b}_{n-1}-n+\beta &=& \frac{an}{\hat{a}_n^2},   \\
  (\hat{a}_{n+1}^2-a)(\hat{a}_n^2-a) &=& a (\hat{b}_n-n) (\hat{b}_n-n+\beta-1),
\end{eqnarray*}
with initial conditions
\[   \hat{a}_0^2=0, \qquad   \hat{b}_0 = \frac{\sqrt{a} I_{-\beta}(2\sqrt{a})}{I_{1-\beta}(2\sqrt{a})}, \]
where $I_\nu$ is the modified Bessel function
\[   I_\nu(z) = \sum_{k=0}^\infty \frac{(z/2)^{2k+\nu}}{k! \Gamma(k+\nu+1)}. \]
\end{theorem}

\begin{proof}
We will give two ways to prove the result. The first way is to use
the relation (\ref{eq:PQ}), which readily gives $\hat{a}_n^2 =
a_n^2$ and $\hat{b}_n = b_n+1-\beta$, where $a_n$ and $b_n$ are the
recurrence coefficients for the polynomials $P_n(x;a,2-\beta)$. We
can therefore use Theorem \ref{thm:1}, but with $\beta$ replaced by
$2-\beta$, and substitute $a_n^2$ by $\hat{a}_n^2$ and $b_n$ by
$\hat{b}_n-1+\beta$. This indeed leaves the equations unchanged. The
initial conditions are $\hat{a}_0^2=0$ and
\[  \hat{b}_0 = b_0 +1-\beta = \frac{\sqrt{a}I_{2-\beta}(2\sqrt{a}) + (1-\beta) I_{1-\beta}(2\sqrt{a})}{I_{1-\beta}(2\sqrt{a})}. \]
The latter expression can be simplified by using a well-known recurrence relation for modified Bessel functions
\cite[Eq.~10.29.1 on p.~251]{NIST}.

An alternative way is to repeat the proof of Theorem \ref{thm:1},
but now on the lattice $\mathbb{N}+1-\beta$. Observe that we only
used the Pearson equation (\ref{eq:Pearson}) and the boundary
condition $w(-1)=0$ in the proof of Theorem \ref{thm:1} to arrive at
the discrete equations (\ref{eq:rec1})--(\ref{eq:rec3}). In the
present case the Pearson equation is still valid and we now use the
boundary condition $w(-\beta)=0$, which allows us to use summation
by parts without the boundary terms. The only difference is that the
moments are now
\begin{eqnarray*}
\hat{m}_0 &=& \sum_{k=0}^\infty w(k+1-\beta) \\
    &=& \Gamma(\beta) a^{1-\beta} \sum_{k=0}^\infty \frac{a^k}{k! \Gamma(k+2-\beta)}
    = \Gamma(\beta) (\sqrt{a})^{1-\beta} I_{1-\beta}(2\sqrt{a}),
\end{eqnarray*}
\begin{eqnarray*}
  \hat{m}_1 &=& \sum_{k=0}^\infty (k+1-\beta) w(k+1-\beta) \\
      &=& \Gamma(\beta) a^{1-\beta} \sum_{k=0}^\infty \frac{a^k}{k! \Gamma(k+1-\beta)}
    = \Gamma(\beta) (\sqrt{a})^{2-\beta} I_{-\beta}(2\sqrt{a}),
\end{eqnarray*}
so that $\hat{b}_0 = \hat{m}_1/\hat{m}_0 = \sqrt{a} I_{-\beta}(2\sqrt{a})/I_{1-\beta}(2\sqrt{a})$.
\end{proof}

\subsection{Combining both lattices}
Now we can combine the two lattices and study orthogonal polynomials on the bi-lattice $\mathbb{N} \cup (\mathbb{N}+1-\beta)$.
We use the orthogonality measure $\mu=c_1 \mu_1 + c_2 \mu_2$, where $c_1,c_2 >0$, $\mu_1$ is the discrete measure on
$\mathbb{N}$ with weights $w_k=w(k)$, and $\mu_2$ is the discrete measure on $\mathbb{N}+1-\beta$ with weights $v_k=w(k+1-\beta)$.
This discrete measure depends on two parameters $c_1,c_2$, but the orthogonal polynomials will only depend on their ratio
$t=c_2/c_1 >0$. Let $0 < \beta < 2$, then $\mu$ is a positive measure and the monic orthogonal polynomials $R_n(x)=R_n(x;a,\beta,t)$
satisfy the three-term recurrence relation
\begin{equation}  \label{eq:Rrec}
    xR_n(x) = R_{n+1}(x) + \tilde{b}_n R_n(x) + \tilde{a}_n^2 R_{n-1}(x).
\end{equation}
It is remarkable that these recurrence coefficients again satisfy the same non-linear recurrence relations (discrete Painlev\'e equations)
as in Theorem \ref{thm:1} and Theorem \ref{thm:2}, but with an initial condition depending on the parameter $t = c_2/c_1$.

\begin{theorem}  \label{thm:3}
The recurrence coefficients for the orthogonal polynomials defined by
\begin{equation}  \label{eq:Rbeta}
  c_1 \sum_{k=0}^\infty R_n(k)R_m(k) w_k + c_2 \sum_{k=0}^\infty R_n(k+1-\beta)R_m(k+1-\beta)v_k = 0, \qquad m \neq n,
\end{equation}
 satisfy the system of equations
\begin{eqnarray*}
   \tilde{b}_n+\tilde{b}_{n-1}-n+\beta &=& \frac{an}{\tilde{a}_n^2},   \\
  (\tilde{a}_{n+1}^2-a)(\tilde{a}_n^2-a) &=& a (\tilde{b}_n-n) (\tilde{b}_n-n+\beta-1),
\end{eqnarray*}
with initial conditions
\[   \tilde{a}_0^2=0, \qquad   \tilde{b}_0 = \sqrt{a} \frac{I_{\beta}(2\sqrt{a})+t I_{-\beta}(2\sqrt{a})}{I_{\beta-1}(2\sqrt{a})+
tI_{1-\beta}(2\sqrt{a})}, \]
where $I_\nu$ is the modified Bessel function
\[   I_\nu(z) = \sum_{k=0}^\infty \frac{(z/2)^{2k+\nu}}{k! \Gamma(k+\nu+1)}. \]
\end{theorem}

\begin{proof}
Going through the proof of Theorem \ref{thm:1} we observe that only
the Pearson equation (\ref{eq:Pearson}) is used, together with
summation by parts. The boundary conditions $w(-1)=0$ and
$w(-\beta)=0$ ensure that this summation by parts does not leave any
boundary terms to evaluate at $-1$ (for the first sum) or $-\beta$
(for the second sum). Hence the recurrence coefficients will satisfy
the same non-linear recurrence relations as in Theorem \ref{thm:1}
and \ref{thm:2}. The only difference is that the initial condition
for $\tilde{b}_0=\tilde{m}_1/\tilde{m}_0$ now is in terms of the
moments
\[    \tilde{m}_0 = c_1 m_0 + c_2 \hat{m_0} = \Gamma(\beta)(\sqrt{a})^{1-\beta} \left( c_1 I_{\beta-1}(2\sqrt{a})+c_2 I_{1-\beta}(2\sqrt{a}) \right), \]
and
\[   \tilde{m}_1 = c_1 m_0 + c_2 \hat{m_0} = \Gamma(\beta)(\sqrt{a})^{2-\beta} \left( c_1 I_{\beta}(2\sqrt{a})+c_2 I_{-\beta}(2\sqrt{a}) \right). \]
\end{proof}

We have now identified special solutions of the discrete system
(\ref{eq:rec1})--(\ref{eq:rec3}) with initial value $a_0^2=0$ and
\[   b_0(t) = \sqrt{a} \frac{I_{\beta}(2\sqrt{a})+t I_{-\beta}(2\sqrt{a})}{I_{\beta-1}(2\sqrt{a})+tI_{1-\beta}(2\sqrt{a})} \]
which depends on one parameter $t \in (0,\infty)$.
If we use the relation \cite[Eq.~10.27.2 on p.~251]{NIST}
\[  I_{-\nu}(z) = I_{\nu}(z) + \frac{2}{\pi} \sin \nu\pi K_{\nu}(z),  \]
where $K_\nu$ is the other modified Bessel function, then the initial condition can also be written as
\[     b_0 = \sqrt{a} \frac{I_\beta(2\sqrt{a}) + s K_\beta(2\sqrt{a})}{I_{\beta-1}(2\sqrt{a}) - s K_{\beta-1}(2\sqrt{a})}, \qquad
  s = \frac{2t}{1+t} \frac{\sin \beta\pi}{\pi}. \]
Observe that
\begin{eqnarray*}
 \frac{I_\beta(2\sqrt{a})}{I_{\beta-1}(2\sqrt{a})} &<& \frac{I_{-\beta}(2\sqrt{a})}{I_{1-\beta}(2\sqrt{a})}, \qquad \textrm{if $0 < \beta < 1$}, \\
 \frac{I_\beta(2\sqrt{a})}{I_{\beta-1}(2\sqrt{a})} &>& \frac{I_{-\beta}(2\sqrt{a})}{I_{1-\beta}(2\sqrt{a})}, \qquad \textrm{if $1 < \beta < 2$}, \\
\end{eqnarray*}
so that $b_0(t)$ is a monotonically increasing function of $t \in [0,\infty]$ when $0 < \beta < 1$ and
a monotonically decreasing function when $1 < \beta < 2$. One can use the Wronskian formula \cite[Eq.~10.28.1 on p.~251]{NIST} and the
product formula \cite[Eq.~10.32.15 on p.~253]{NIST} to verify this. For each initial value in $[b_0(0),b_0(\infty)]$ (when $0 < \beta < 1$)
or $[b_0(\infty),b_0(0)]$ (when $1 < \beta < 2$) the solution therefore corresponds to recurrence coefficients of orthogonal polynomials
with a positive measure on the real line (in fact on the bi-lattice $\mathbb{N} \cup (\mathbb{N}+1-\beta)$) whenever $a >0$, and
hence this solution satisfies certain positivity constraints:  $a_n^2 > 0$ for $n \geq 1$ and $b_n > \min(0,1-\beta)$
for all $n \geq 0$. Moreover, since the
orthogonal polynomials are on  a discrete set, we can use the familiar fact for orthogonal polynomials that between two zeros there is at least one
point of the bi-lattice \cite[Theorem 4.1 on p.~59]{Chihara}. The sum of the zeros $x_{1,n} < x_{2,n} < \cdots < x_{n,n}$ of $R_n$ hence satisfies
\[   \sum_{k=1}^n x_{k,n}  > \sum_{k=0}^{n-1} y_k  \]
where $y_0 < y_1 < y_2 < \cdots$ are the points in the bi-lattice $\mathbb{N} \cup (\mathbb{N}+1-\beta)$, i.e.,
\[  \begin{cases}
   y_{2k} = k, \ y_{2k+1} = k+1-\beta, & \textrm{if $0 < \beta < 1$}, \\
   y_{2k} = k+1-\beta, \ y_{2k+1} = k, & \textrm{if $1 < \beta < 2$}, \\
   y_k = k, & \textrm{if $\beta=1$}.
    \end{cases}  \]
The sum of the zeros is the trace of the truncated Jacobi matrix containing the recurrence coefficients, hence
\[   \sum_{k=1}^n x_{k,n} = \sum_{k=0}^{n-1} b_k , \]
so that we get the constraint
\begin{eqnarray*}
   \sum_{k=0}^{2n-1} b_k &>& n(n-1)+n(1-\beta),  \qquad \textrm{if $0 < \beta < 2$ and $\beta \neq 1$}, \\
    \sum_{k=0}^{n-1} b_k &>& n(n-1)/2,  \qquad \textrm{if $\beta=1$}.
\end{eqnarray*}
The sum $\sum_{k=0}^{n-1} b_k$ behaves like $n^2/4+ \mathcal{O}(n)$
as $n \to \infty$ when $\beta \neq 1$ and as $n^2/2 +
\mathcal{O}(n)$ when $\beta=1$. This result follows from
\cite[Theorem 2.2]{KuijlWVA} (see in particular the example of
Charlier polynomials on p.~200). This difference in the behavior of
the recurrence coefficients is easily explained since for $\beta=1$
the two lattices coincide and we just have one lattice $\mathbb{N}$.
A similar situation appears when $t \to 0$ (which gives the
orthogonal polynomials $P_n$ on the lattice $\mathbb{N}$) and when
$t \to \infty$, which gives the polynomials $Q_n$ on the lattice
$\mathbb{N}+1-\beta$. We believe that these extreme cases ($t = 0$,
$t = \infty$ and $\beta=1$) correspond to the only solution of
(\ref{eq:rec1})--(\ref{eq:rec3}) with $a_0^2=0$ for which $a_n^2 >0$
for all $n >0$ and $b_n = n + \mathcal{O}(1)$ as $n \to \infty$.

Another special case occurs when $\beta=1/2$. In this case the bi-lattice is equally spaced and equal to $\frac12 \mathbb{N}$. The lattice points
are $y_k = k/2$ and the weight $w$ at these points is
\[     w(k/2) = \frac{\sqrt{\pi} a^{k/2}}{\Gamma((k+1)/2)\Gamma(k/2+1)} = \frac{2^k a^{k/2}}{\Gamma(k+1)} \]
where we used Legendre's duplication formula \cite[Eq.~5.5.5 on p.~138]{NIST} for the gamma function. This means that if we take $c_1=c_2$, then
the weights of the measure $\mu=\mu_1+\mu_2$ are precisely the Charlier weights (Poisson distribution) with parameter $2\sqrt{a}$.
The orthogonal polynomials are therefore Charlier polynomials \cite[Chapter VI, \S 1]{Chihara} but with a scaling:
\[     R_n(x;a,1/2,1) = 2^{-n} \hat{C}_n(2x;2\sqrt{a}), \]
where $\hat{C}_n(x;a)$ are the monic Charlier polynomials with parameter $a$. The recurrence coefficients of Charlier polynomials
$\hat{C}_n(x;a)$ are known to be $a_n^2 = na$ and $b_n = n+a$, so for the scaled polynomials $2^{-n}\hat{C}_n(2x;a)$ they are
$a_n^2 = na/4$ and $b_n=(n+a)/2$, and for our polynomials $R_n(x) = 2^{-n} \hat{C}_n(2x;2\sqrt{a})$ we thus have
\begin{equation}  \label{eq:beta1/2}
    \tilde{a}_n^2 = \frac{n\sqrt{a}}{2}, \quad \tilde{b}_n = \frac{n}{2} + \sqrt{a}.
\end{equation}
One can indeed verify that (\ref{eq:beta1/2}) is a solution of
(\ref{eq:rec1})--(\ref{eq:rec3}) with initial condition
\[   \tilde{b}_0 = \sqrt{a} \frac{I_{1/2}(2\sqrt{a})+I_{-1/2}(2\sqrt{a})}{I_{-1/2}(2\sqrt{a})+I_{1/2}(2\sqrt{a})} = \sqrt{a}. \]
So there are special values of $\beta$ for which the non-linear equations have simple solutions (as a function of $n$).

\section{A generalized Meixner weight}
\subsection{The lattice $\mathbb{N}$}
We consider the sequence $(p_n)_n$ of polynomials, \textit{orthonormal} with
respect to the weights
\begin{equation}
w_k=\frac{(\gamma)_ka^k}{k!(\beta)_k}\label{GenMeixWeight}, \qquad k
\in \mathbb{N} = \{0,1,2,3,\ldots \}.
\end{equation}
These polynomials satisfy
\[\sum_{k=0}^\infty p_n(k)p_m(k)w_k=\delta_{n,m}, \]
where $\delta_{n,m}$ is the Kronecker delta.  These weights are
positive whenever $\beta, \gamma, a >0$. Notice that the special
case $\beta=1$ was studied in \cite{LiesWva}, where the authors
showed that the recurrence coefficients $a_n,b_n$ satisfy a limiting
case of the asymmetric discrete Painlev\'e equation
$\alpha$-dP$_{\rm IV}$.  The special case $\beta=\gamma$ gives the
well-known Charlier weight. We can write these weights in terms of a
weight function $w$ given by
\begin{equation}  \label{eq:Meixw}
     w(x) = \frac{\Gamma(\beta)}{\Gamma(\gamma)} \frac{\Gamma(\gamma+x) a^x}{\Gamma(\beta+x) \Gamma(x+1)}
\end{equation}
so that $w_k = w(k)$ for $k \in \mathbb{N}$.

We will prove the following result on the recurrence coefficients of the sequence of orthonormal
polynomials.
\begin{theorem} \label{GenMeixThmOriginalLattice}
The orthonormal polynomials with respect to the
weights \eqref{GenMeixWeight} on the lattice $\mathbb{N}$, with $\gamma,\beta,a>0$,  
satisfy the three-term recurrence relation
\[   xp_n(x)=a_{n+1}p_{n+1}(x)+b_np_n(x)+a_np_{n-1}(x),\]
where the recurrence coefficients are given by the initial conditions
\[   a_0=0, \qquad
     b_0=\frac{\gamma a}{\beta} \frac{M(\gamma+1,\beta+1,a)}{M(\gamma,\beta,a)},\]
where $M(a,b,z)$ is the confluent hypergeometric function
\[    M(a,b,z) = \sum_{k=0}^\infty \frac{(a)_k}{(b)_k k!} z^k = {}_1F_1(a;b;z), \]
and $a_n^2= na-(\gamma-1)u_n$, $b_n =
n+\gamma-\beta+a-(\gamma-1)v_n/a$, where $(u_n,v_n)_{n \in
\mathbb{N}}$ satisfy the system of non-linear equations
\begin{eqnarray}
     (u_n+v_n)(u_{n+1}+v_n)&=&\frac{\gamma-1}{a^2}v_n(v_n-a)\left(v_n-a\frac{\gamma-\beta}{\gamma-1}\right),  \label{eq:uvrec1}  \\
     (u_n+v_n)(u_n+v_{n-1})&=&\frac{u_n}{u_n-\frac{an}{\gamma-1}}(u_n+a)\left(u_n+a\frac{\gamma-\beta}{\gamma-1}\right).
                        \label{eq:uvrec2}
\end{eqnarray}
\end{theorem}

\begin{proof}
We will use the technique of ladder operators, proposed by Chen and
Ismail \cite{ChenIsmail}. The ladder operators are defined by
\begin{eqnarray*}
A_n(x)&=&a_n\sum_{\ell=0}^\infty p_n(\ell)p_n(\ell-1)\frac{u(x+1)-u(\ell)}{x+1-\ell}w(\ell),\\
B_n(x)&=&a_n\sum_{\ell=0}^\infty
p_n(\ell)p_{n-1}(\ell-1)\frac{u(x+1)-u(\ell)}{x+1-\ell}w(\ell),
\end{eqnarray*}
where $u(x)=-1+\frac{w(x-1)}{w(x)}$.  These ladder operators satisfy
\begin{equation}\label{GenMeixLadderOperatorsOrLatt}
A_n(x)p_{n-1}(x)-B_n(x)p_n(x)=p_n(x+1)-p_n(x).
\end{equation}
This can be seen by simplifying the left hand side using the
definitions for $A_n$, $B_n$ and $u$, and the Christoffel-Darboux
identity
\[\sum_{j=0}^np_j(x)p_j(y)=a_{n+1}\frac{p_{n+1}(x)p_n(y)-p_n(x)p_{n+1}(y)}{x-y}.\] After doing
this we obtain
\[A_n(x)p_{n-1}(x)-B_n(x)p_n(x)=\sum_{j=0}^{n-1}p_j(x)\sum_{\ell=0}^\infty
p_n(\ell+1)p_j(\ell)w(\ell).\]  Exactly the same expression is found
by writing $p_n(x+1)-p_n(x)$ as a Fourier series:
\[p_n(x+1)-p_n(x)=\sum_{j=0}^{n-1}\alpha_{n,j}p_j(x).\]
Equation \reff{GenMeixLadderOperatorsOrLatt} gives rise to the
compatibility relations
\begin{equation}\label{CompRel1}B_n(x)+B_{n+1}(x)=\frac{x-b_n}{a_n}A_n(x)-u(x+1)+\sum_{j=0}^n\frac{A_j(x)}{a_j}\end{equation}
and
\begin{equation}\label{CompRel2}a_{n+1}A_{n+1}(x)-a_n^2\frac{A_{n-1}(x)}{a_{n-1}}=(x-b_n)B_{n+1}(x)-(x+1-b_n)B_n(x)+1.\end{equation}
An easy calculation gives that
\[\frac{u(x+1)-u(\ell)}{x+1-\ell}=\frac{1}{a(\gamma+\ell-1)}\left(\ell+\frac{(\gamma-1)(x+\beta)}{\gamma+x}\right).\]
Inserting this into the ladder operators, we find that
\begin{equation}\label{GenMeixLadderAn}A_n(x)=\frac{a_n}{a}R_n+\frac{a_n}{a}\frac{x+\beta}{x+\gamma}T_n\end{equation}
and
\begin{equation}\label{GenMeixLadderBn}B_n(x)=\frac{1}{a}r_n+\frac{1}{a}\frac{x+\beta}{x+\gamma}t_n,\end{equation}
where
\[R_n=\sum_{\ell=0}^\infty p_n(\ell)p_n(\ell-1)\frac{\ell}{\gamma+\ell-1}w(\ell),\]
\[T_n=\sum_{\ell=0}^\infty p_n(\ell)p_n(\ell-1)\frac{\gamma-1}{\gamma+\ell-1}w(\ell),\]
\[r_n=a_n\sum_{\ell=0}^\infty p_n(\ell)p_{n-1}(\ell-1)\frac{\ell}{\gamma+\ell-1}w(\ell),\]
\[t_n=a_n\sum_{\ell=0}^\infty p_n(\ell)p_{n-1}(\ell-1)\frac{\gamma-1}{\gamma+\ell-1}w(\ell).\]
We use the expressions
\eqref{GenMeixLadderAn}--\eqref{GenMeixLadderBn} in the compatibility
relations, and compare the coefficients of $x^2,x$ and $1$.  In
this way we obtain a system of six equations:
\begin{equation}\label{GenMeix1of6}
 0=R_n+T_n-1,
\end{equation}
\begin{equation}\label{GenMeix2of6}
r_n+t_n+r_{n+1}+t_{n+1}=\gamma R_n-b_nR_n-b_nT_n+\beta
T_n+a-\beta-1+\sum_{j=0}^nR_j+\sum_{j=0}^nT_j,
\end{equation}
\begin{equation}\label{GenMeix3of6}
\gamma r_n+\beta t_n+\gamma r_{n+1}+\beta t_{n+1}=-b_n\gamma
R_n-b_n\beta T_n+\gamma a-\beta+\gamma\sum_{j=0}^nR_j+
\beta\sum_{j=0}^nT_j,
\end{equation}
\begin{equation}\label{GenMeix4of6}
0=r_{n+1}+t_{n+1}-r_n-t_n,
\end{equation}
\begin{multline}\label{GenMeix5of6}
a_{n+1}^2R_{n+1}+a_{n+1}^2T_{n+1}-a_n^2R_{n-1}-a_n^2T_{n-1}\\
=(\gamma-b_n) r_{n+1}+(\beta-b_n) t_{n+1}
-(1-b_n+\gamma)r_n-(1-b_n+\beta)t_n+a,
\end{multline}
\begin{multline}\label{GenMeix6of6}
\gamma a_{n+1}^2R_{n+1}+\beta a_{n+1}^2T_{n+1}-\gamma a_n^2R_{n-1}-\beta a_n^2T_{n-1}\\
=-\gamma b_nr_{n+1}-\beta
b_nt_{n+1}-\gamma(1-b_n)r_n-\beta(1-b_n)t_n+\gamma a.
\end{multline}
By \eqref{GenMeix1of6} we can substitute $R_n$ by $1-T_n$, and by
\eqref{GenMeix4of6} (and since $r_0=t_0=0$ by definition) we can
substitute $r_n$ by $-t_n$.  If we apply these substitutions,
\eqref{GenMeix2of6} gives an expression for $b_n$ as a function of
$T_n$:
\begin{equation}\label{GenMeixbn}b_n=\gamma-(\gamma-\beta)T_n+a+n-\beta,\end{equation} and
\eqref{GenMeix3of6} gives
\[(\gamma-\beta)(t_n+t_{n+1})=\gamma
b_n-(\gamma-\beta)b_nT_n-\gamma
a+\beta-\gamma(n+1)+(\gamma-\beta)\sum_{j=0}^nT_j.\] Substituting
the newly found expression for $b_n$ in this equation, we obtain
\begin{equation}\label{GenMeixtntn+1}t_n+t_{n+1}=(\gamma-\beta)T_n^2-(2\gamma+a+n-\beta-1)T_n+\gamma-1+\sum_{j=0}^{n-1}T_j.\end{equation}
On the other hand, elimination of $R_n$ and $r_n$ in
\eqref{GenMeix6of6} gives
\begin{equation}\label{GenMeixanTnbntn}\gamma(a_{n+1}^2-a_n^2)-(\gamma-\beta)a_{n+1}^2T_{n+1}+(\gamma-\beta)a_n^2T_{n-1}=(\gamma-\beta)b_nt_{n+1}+
(\gamma-\beta)(1-b_n)t_n+\gamma a\end{equation} and in
\eqref{GenMeix5of6} it gives
\begin{equation}\label{GenMeixanan+1tntn+1}a_{n+1}^2-a_n^2=(\gamma-\beta)(t_n-t_{n+1})+a,\end{equation}
which after taking a telescoping sum gives us an expression for the
recurrence coefficients $a_n$ as a function of
$t_n$:\begin{equation}\label{GenMeixan}a_n^2=na-(\gamma-\beta)t_n.\end{equation}
Notice that since $a_n\geq 0$ for all $n$, the knowledge of $t_n$
will imply the knowledge of $a_n$.  Using
\eqref{GenMeixanan+1tntn+1} and the expression for $b_n$ in
\eqref{GenMeixanTnbntn}, we obtain
\begin{equation}\label{GenMeixan2Tn-1}a_n^2T_{n-1}-a_{n+1}^2T_{n+1}=(b_n+\gamma)t_{n+1}+(1-b_n-\gamma)t_n,\end{equation}
in which substitution of \eqref{GenMeixan} gives us another relation
between terms of the sequences $(t_n)_n$ and $(T_n)_n$:
\begin{eqnarray*}-a(n+1)T_{n+1}+anT_{n-1}&=&t_{n+1}(n+a+2\gamma-\beta)-t_n(n+a+2\gamma-\beta-1)\\
&&-t_{n+1}(\gamma-\beta)(T_n+T_{n+1})+t_n(\gamma-\beta)(T_n+T_{n-1}).\end{eqnarray*}
Taking a telescopic sum, we get
\begin{equation}\label{GenMeixTn+Tn-1}(T_n+T_{n-1})(an-(\gamma-\beta)t_n)=-(a+2\gamma+n-\beta-1)t_n+a\sum_{j=0}^{n-1}T_j.\end{equation}
On the other hand, if we multiply \eqref{GenMeixan2Tn-1} by $T_n$
and we use \eqref{GenMeixbn} and \eqref{GenMeixtntn+1}, we obtain
\[a_{n+1}^2T_{n+1}T_n-a_n^2T_nT_{n-1}=t_{n+1}^2-t_n^2-(\gamma-1)(t_{n+1}-t_n)-t_{n+1}\sum_{j=0}^nT_j+t_n\sum_{j=0}^{n-1}T_j,\]
which by taking a telescoping sum gives
\begin{equation}\label{GenMeixan2TnTn-1}a_n^2T_nT_{n-1}=t_n\left(t_n-\gamma+1-\sum_{j=0}^{n-1}T_j\right).\end{equation}
If we multiply \eqref{GenMeixTn+Tn-1} by $T_n$ and
\eqref{GenMeixtntn+1} by $t_n$, comparison of these two expressions
gives \[anT_n^2-t_nt_{n+1}=aT_n\sum_{j=0}^{n-1}T_j.\]  Using in this
equality the expression for $\sum_{j=0}^{n-1}T_j$ from
\eqref{GenMeixtntn+1}, we obtain a relation between $T_n$, $t_n$ and
$t_{n+1}$:
\begin{equation}\label{GenMeixtnTnsystem1}aT_n(T_n-1)(T_n(\gamma-\beta)-\gamma+1)=(t_n+aT_n)(t_{n+1}+aT_n).\end{equation}
Finally, multiplying \eqref{GenMeixTn+Tn-1} by $at_n$, and using
once again \eqref{GenMeixan2TnTn-1} to eliminate the sum of $T_j$,
and \eqref{GenMeixan} to eliminate $a_n^2$, we obtain a relation
between $T_n$, $T_{n-1}$ and
$t_n$:\begin{equation}\label{GenMeixtnTnsystem2}(aT_n+t_n)(aT_{n-1}+t_n)=\frac{t_n(a^2(\gamma-1)+a(2\gamma-\beta-1)t_n+(\gamma-\beta)t_n^2)}
{(\gamma-\beta)t_n-an}.\end{equation} The substitution
\[u_n=\frac{\gamma-\beta}{\gamma-1}t_n,\qquad
v_n=\frac{\gamma-\beta}{\gamma-1}aT_n, \] allows us to write the
system of two equations in a more symmetrical form:
\begin{eqnarray*}
(u_n+v_n)(u_{n+1}+v_n)&=&\frac{\gamma-1}{a^2}v_n(v_n-a)\left(v_n-a\frac{\gamma-\beta}{\gamma-1}\right), \\
(u_n+v_n)(u_n+v_{n-1})&=&\frac{u_n}{u_n-\frac{an}{\gamma-1}}(u_n+a)\left(u_n+a\frac{\gamma-\beta}{\gamma
-1}\right),
\end{eqnarray*}
which is the system (\ref{eq:uvrec1})--(\ref{eq:uvrec2}) given in
the Theorem. The initial conditions are given by $u_0=0$, since
$t_0=0$, and
$v_0=\frac{a}{\gamma-1}(\gamma-\beta+a-\frac{m_1}{m_0})$ since $T_0$
and $b_0$ are related by \eqref{GenMeixbn}, and since $b_0=m_1/m_0$,
where $m_k$ is the $k$-th moment of the weights $(w_k)_{k \in
\mathbb{N}}$.  The first two moments are easily recognized as
confluent hypergeometric functions:
\[ m_0=\sum_{k=0}^\infty w_k= M(\gamma,\beta,a)  \]
and
\[ m_1=\sum_{k=0}^\infty  kw_k  = \frac{\gamma a}{\beta}  M(\gamma+1,\beta+1,a). \]
\end{proof}

The system (\ref{eq:uvrec1})--(\ref{eq:uvrec2}) can be seen as a limiting case of the asymmetric discrete
Painlev\'e-IV equation $\alpha$-dP$_{\rm IV}$ (\cite[$E_6^{\delta}$ on p.~296]{GR}, \cite[Appendix A.3]{WVA}) given by
\begin{eqnarray*}
           \ds (X_n+Y_n)(X_{n+1}+Y_n)&=&\frac{(Y_n-A)(Y_n-B)(Y_n-C)(Y_n-D)}{(Y_n+\Gamma-Z_n) (Y_n-\Gamma-Z_n)},\\
           \ds (X_n+Y_n)(X_n+Y_{n-1})&=&\frac{(X_n+A)(X_n+B)(X_n+C)(X_n+D)}{(X_n+\Delta-Z_{n+1/2}) (X_n-\Delta-Z_{n+1/2})},
\end{eqnarray*}
where $A+B+C+D=0$ and $Z_n$ is a polynomial of degree one in $n$.
Indeed, let
\[X_n=u_n-\frac{1}{\epsilon},\ Y_n=v_n+\frac{1}{\epsilon},\
Z_n=\frac{a}{\gamma-1}\left(n-\frac{1}{2}\right)+\frac{1}{\epsilon},\]
and
\[A=\frac{1}{\epsilon},\
B=-\frac{3}{\epsilon}-a-a\frac{\gamma-\beta}{\gamma-1},\
C=a+\frac{1}{\epsilon},\
D=\frac{1}{\epsilon}+a\frac{\gamma-\beta}{\gamma-1},\]
\[\Gamma^2=\frac{-4a^2}{(\gamma-1)\epsilon}, \quad
\Delta=\frac{2}{\epsilon},\]
then letting $\epsilon$ tend to zero
gives the system (\ref{eq:uvrec1})--(\ref{eq:uvrec2}).


\subsection{The shifted lattice $\mathbb{N}+1-\beta$}
The weight function $w$ has zeros at the negative
integers $-1,-2,\ldots$, but also at $-\beta,-1-\beta,\ldots$. Therefore it
makes sense to consider the weight function on the shifted lattice
$\mathbb{N}+1-\beta$.
It is easy to verify that the weight on this shifted lattice is, up
to a constant factor, equal to the weight on the original lattice
$\mathbb{N}$, with different parameters. Indeed, if we denote
\[ w_{\gamma,\beta,a}(x)=\frac{\Gamma(\beta)}{\Gamma(\gamma)} \frac{\Gamma(\gamma+x)a^x}{\Gamma(x+1)\Gamma(\beta+x)}, \]
then
\[w_{\gamma,\beta,a}(k+1-\beta)=a^{1-\beta}\frac{\Gamma(\beta)\Gamma(\gamma+1-\beta)}
{\Gamma(2-\beta)\Gamma(\gamma)}w_{\gamma+1-\beta,2-\beta,a}(k).\]
Hence the corresponding orthonormal polynomials $q_n$, which satisfy
\begin{equation}  \label{eq:qortho}
\sum_{k=0}^\infty q_n(k+1-\beta)q_m(k+1-\beta)w(k+1-\beta) = \delta_{n,m} ,
\end{equation}
are, up to a constant factor, equal to the polynomials $p_n$,
shifted in both the variable $x$ and the parameters $\gamma$ and
$\beta$:
\begin{equation}  \label{eq:pq}
q_n^{\gamma,\beta,a}(x)=\zeta_np_n^{\gamma+1-\beta,2-\beta,a}(x+\beta-1).
\end{equation}
The weights $(w(k+1-\beta))_{k \in \mathbb{N}}$ are positive when $a
> 0$, $\gamma+1-\beta >0$ and $2-\beta>0$, i.e., whenever $a>0$,
$\beta < 2$ and $\gamma > \beta-1$. Remarkably, the recurrence
coefficients in the recurrence relation
\[xq_n(x)=\hat{a}_{n+1}q_{n+1}(x)+\hat{b}_nq_n(x)+\hat{a}_nq_{n-1}(x)\]
satisfy exactly the same system of difference equations as the
polynomials $p_n$, only the initial condition for $\hat{b}_0$
changes.
\begin{theorem} The recurrence coefficients for the orthonormal
polynomials $q_n$, defined by (\ref{eq:qortho}) are given by
$\hat{a}_n^2 = na - (\gamma-1)\hat{u}_n$, $\hat{b}_n =
n+\gamma-\beta+a - (\gamma-1)\hat{v}_n/a$, where
$(\hat{u}_n,\hat{v}_n)_{n \in \mathbb{N}}$ satisfy the system of
equations
\begin{eqnarray*}
         (\hat{u}_n+\hat{v}_n)(\hat{u}_{n+1}+\hat{v}_n)&=&\frac{\gamma-1}{a^2}\hat{v}_n(\hat{v}_n-a)
                          \left(\hat{v}_n-a\frac{\gamma-\beta}{\gamma-1}\right) \\
         (\hat{u}_n+\hat{v}_n)(\hat{u}_n+\hat{v}_{n-1})&=&\frac{\hat{u}_n}{\hat{u}_n-\frac{an}{\gamma-1}}(\hat{u}_n+a)
                          \left(\hat{u}_n+a\frac{\gamma-\beta}{\gamma-1}\right),
\end{eqnarray*}
with initial conditions
\[  \hat{a}_0=0, \quad
    \hat{b}_0= (1-\beta) \frac{M(\gamma-\beta+1,1-\beta,a)}{M(\gamma-\beta+1,2-\beta,a)}, \]
where $M(a,b,z)$ is the confluent hypergeometric function.
\end{theorem}
\begin{proof}
The proof of Theorem~\ref{GenMeixThmOriginalLattice} can easily be
adapted to this case, where the ladder operators are now given by
\begin{eqnarray*}
\hat{A}_n(x)&=&\hat{a}_n\sum_{\ell=0}^\infty q_n(\ell+1-\beta)q_n(\ell-\beta)\frac{u(x+1)-u(\ell+1-\beta)}{x-\ell+\beta}w(\ell+1-\beta),\\
\hat{B}_n(x)&=&\hat{a}_n\sum_{\ell=0}^\infty
q_n(\ell+1-\beta)q_{n-1}(\ell-\beta)\frac{u(x+1)-u(\ell+1-\beta)}{x-\ell+\beta}w(\ell+1-\beta).
\end{eqnarray*}  They satisfy
\[\hat{A}_n(x)q_{n-1}(x)-\hat{B}_n(x)q_n(x)=q_n(x+1)-q_n(x),\]
leading to exactly the same compatibility relations for $\hat{A}_n$
and $\hat{B}_n$ as in \reff{CompRel1}--\reff{CompRel2}.  Hence we
obtain the same difference relation for the recurrence coefficients
$\hat{a}_n$ and $\hat{b}_n$ as in the case of the lattice
$\mathbb{N}$.  The only difference is the initial condition for
$\hat{b}_0$, which can again be obtained by calculating
$\hat{m}_1/\hat{m}_0$.

Alternatively, if we use the relation (\ref{eq:pq}), then we see that
\[   \hat{a}_n^2 = a_n^2(\gamma-\beta+1,2-\beta,a), \quad \hat{b}_n = b_n(\gamma-\beta+1,2-\beta,a) + 1-\beta, \]
where $a_n$ and $b_n$ are the recurrence coefficients of the polynomials in Theorem \ref{GenMeixThmOriginalLattice}
but with $\gamma$ replaced by $\gamma-\beta+1$ and $\beta$ replaced by $2-\beta$. This means that
\[  \hat{u}_n= \frac{\gamma-\beta}{\gamma-1} u_n(\gamma-\beta+1,2-\beta,a), \quad
    \hat{v}_n = \frac{\gamma-\beta}{\gamma-1} v_n(\gamma-\beta+1,2-\beta,a) \]
and if we insert this in (\ref{eq:uvrec1})--(\ref{eq:uvrec2}), then we retrieve the same non-linear system
of equations for $(\hat{u}_n,\hat{v}_n)_{n \in \mathbb{N}}$. The initial condition is
\[ \hat{b}_0 = b_0(\gamma-\beta+1,2-\beta,a) + 1-\beta, \]
which becomes
\[   \hat{b}_0 = \frac{a(\gamma-\beta+1)M(\gamma+2-\beta,3-\beta,a)+(2-\beta)(1-\beta)M(\gamma+1-\beta,2-\beta,a)}
                         {(2-\beta)M(\gamma-\beta+1,2-\beta,a)}.  \]
This can be simplified by using some recurrence relations for confluent hypergeometric functions, in particular 13.3.3 and 13.3.4 in
\cite[p.~325]{NIST}, which leads to
\[ \hat{b}_0 = (1-\beta) \frac{M(\gamma-\beta+1,1-\beta,a)}{M(\gamma-\beta+1,2-\beta,a)}  . \]
\end{proof}

\subsection{Combining both lattices}
We can now consider a combination of both lattices by considering
the measure $\mu=\mu_1+t\mu_2$ on the bi-lattice
$\mathbb{N}\cup (\mathbb{N}+1-\beta)$, where $\mu_1$ is the discrete
measure on $\mathbb{N}$ with weights $(w(k))_{k\in \mathbb{N}}$, and $\mu_2$ is the
discrete measure on $\mathbb{N}+1-\beta$ with weights $(w(k+1-\beta))_{k\in \mathbb{N}}$.
In order to have two positive measures $\mu_1$ and $\mu_2$ we impose the conditions
\[   a > 0,\quad 0 < \beta < 2, \quad \gamma > \max(0,\beta-1), \quad t \in [0,\infty].  \]
The orthonormal polynomials $(r_n)_{n \in \mathbb{N}}$ for this measure
depend on $t > 0$.  Once again the recurrence coefficients which
appear in the recurrence relation
\[xr_n(x)=\tilde{a}_{n+1}r_{n+1}(x)+\tilde{b}_nr_n(x)+\tilde{a}_nr_{n-1}(x)\]
satisfy the same system of non-linear equations.
\begin{theorem}
The recurrence coefficients for the orthonormal
polynomials $r_n$, defined by
\[\sum_{k=0}^\infty
r_n(k)r_m(k)w(k)+t\sum_{k=0}^\infty
r_n(k+1-\beta)r_m(k+1-\beta)w(k+1-\beta)=\delta_{n,m} \] are given
by $\tilde{a}_n^2 = na - (\gamma-1)\tilde{u}_n$ and $\tilde{b}_n =
n+\gamma-\beta+a-(\gamma-1)\tilde{v}_n/a$, where
$(\tilde{u}_n,\tilde{v}_n)_{n \in \mathbb{N}}$ satisfy the system of
non-linear equations
\begin{eqnarray*}
         (\tilde{u}_n+\tilde{v}_n)(\tilde{u}_{n+1}+\tilde{v}_n)&=&\frac{\gamma-1}{a^2}\tilde{v}_n(\tilde{v}_n-a)
                          \left(\tilde{v}_n-a\frac{\gamma-\beta}{\gamma-1}\right), \\
         (\tilde{u}_n+\tilde{v}_n)(\tilde{u}_n+\tilde{v}_{n-1})&=&\frac{\tilde{u}_n}{\tilde{u}_n-\frac{an}{\gamma-1}}(\tilde{u}_n+a)
                          \left(\tilde{u}_n+a\frac{\gamma-\beta}{\gamma-1}\right),
\end{eqnarray*}
with initial conditions $\tilde{a}_0=0$ and
\begin{equation}  \label{GenMeixInitCondTilde}
    \tilde{b}_0= \frac{m_1+t\hat{m}_1}{m_0+t\hat{m}_0},
\end{equation}
where
\[   m_0 = M(\gamma,\beta,a), \quad m_1 = \frac{\gamma a}{\beta} M(\gamma+1,\beta+1,a), \]
\[   \hat{m}_0 = \frac{\Gamma(\beta)\Gamma(\gamma-\beta+1)}{\Gamma(\gamma)\Gamma(2-\beta)} a^{1-\beta} M(\gamma-\beta+1,2-\beta,a), \]
\[    \hat{m}_1 = \frac{\Gamma(\beta)\Gamma(\gamma-\beta+1)}{\Gamma(\gamma)\Gamma(1-\beta)} a^{1-\beta} M(\gamma-\beta+1,1-\beta,a). \]
\end{theorem}

\begin{proof}
Once again, the proof of Theorem~\ref{GenMeixThmOriginalLattice} can
be adapted to this case.  Now the ladder operators are given by
\begin{eqnarray*}
\tilde{A}_n(x)&=&\tilde{a}_n\left[\sum_{\ell=0}^\infty
r_n(\ell)r_n(\ell-1)\frac{u(x+1)-u(\ell)}{x+1-\ell}w(\ell)\right.\\
&&\left.+t\sum_{\ell=0}^\infty r_n(\ell+1-\beta)r_n(\ell-\beta)\frac{u(x+1)-u(\ell+1-\beta)}{x-\ell+\beta}w(\ell+1-\beta)\right],\\
\tilde{B}_n(x)&=&\tilde{a}_n\left[\sum_{\ell=0}^\infty
r_n(\ell)r_{n-1}(\ell-1)\frac{u(x+1)-u(\ell)}{x+1-\ell}w(\ell)\right.\\
&&\left.+t\sum_{\ell=0}^\infty
r_n(\ell+1-\beta)r_{n-1}(\ell-\beta)\frac{u(x+1)-u(\ell+1-\beta)}{x-\ell+\beta}w(\ell+1-\beta)\right].
\end{eqnarray*}  They satisfy
\[\tilde{A}_n(x)r_{n-1}(x)-\tilde{B}_n(x)r_n(x)=r_n(x+1)-r_n(x),\]
leading to exactly the same compatibility relations for
$\tilde{A}_n$ and $\tilde{B}_n$ as in
\reff{CompRel1}--\reff{CompRel2}. Hence we obtain the same difference
relation for the recurrence coefficients $\tilde{a}_n$ and
$\tilde{b}_n$ as in the case of the lattice $\mathbb{N}$ or the
shifted lattice $\mathbb{N}+1-\beta$. The only difference is the
initial condition for $\tilde{b}_0$, which can again be obtained by
calculating $\tilde{m}_1/\tilde{m}_0$, where $\tilde{m}_k = m_k + t\hat{m}_k$,
with $m_k$ and $\hat{m}_k$ the $k$-th moment of the measures $\mu_1$ and $\mu_2$
respectively.
\end{proof}

Notice that letting $t$ tend to $0$ or $\infty$, we obtain the
initial condition for $b_0$ of Theorem 3.1 and $\hat{b}_0$ of Theorem 3.2 respectively, which was to be
expected. If we use the formula \cite[Eq.~13.2.42 on p.~325]{NIST}
\[   U(a,b,z) = \frac{\Gamma(1-b)}{\Gamma(a-b+1)} M(a,b,z) + \frac{\Gamma(b-1)}{\Gamma(a)} z^{1-b} M(a-b+1,2-b,z) \]
then the initial condition can also be written as
\[  \tilde{b}_0 = \frac{\gamma a}{\beta} \frac{M(\gamma+1,\beta+1,a) - s \Gamma(\gamma-\beta+1)/\Gamma(-\beta) U(\gamma+1,\beta+1,a)}
             {M(\gamma,\beta,a)-s \Gamma(\gamma-\beta+1)/\Gamma(1-\beta) U(\gamma,\beta,a)}, \qquad s = \frac{t}{1+t},  \]
where $U(a,b,z)$ is the second solution of the confluent hypergeometric differential equation.

Now we have identified a class of solutions to
\reff{eq:uvrec1}--\reff{eq:uvrec2}, with initial conditions
$\tilde{a}_0=0$ and \reff{GenMeixInitCondTilde}, depending on a
parameter $t$.  This initial condition is a homographic function of
$t$.  Using the Wronskian formula \cite[Eq.~13.2.33]{NIST} and the
contiguous relations \cite[Eq.~13.3.15]{NIST} and
\[\frac{az}{b}M(a+1,b+1,z)=(1-b)(M(a,b,z)-M(a,b-1,z))\] from \cite[Eqs.~13.3.3 and 13.3.4]{NIST}, one can prove that
$\tilde{b}_0$ is an increasing function of $t \in \mathbb{R}^+$ when
$0<\beta<1$, and a decreasing function of $t \in \mathbb{R}^+$ when $1 < \beta < 2$.
Hence for each initial value in
$[\tilde{b}_0(0),\tilde{b}_0(\infty)]$ or
$[\tilde{b}_0(\infty),\tilde{b}_0(0)]$, the solution to this
non-linear system corresponds to recurrence coefficients of
orthogonal polynomials on the bi-lattice
$\mathbb{N}\cup(\mathbb{N}+1-\beta)$.

Exactly the same results concerning the asymptotic behavior of the
sum $\sum_{k=0}^{n-1}b_k$ hold as in the case of the generalized
Charlier polynomials: the sum behaves as $n^2/4+\mathcal{O}(n)$ as
$n\rightarrow\infty$ when $\beta\neq 1$, and as
$n^2/2+\mathcal{O}(n)$ when $\beta=1$ or for the limiting cases $t=0$ and $t=\infty$, which correspond
to the latiice $\mathbb{N}$ and $\mathbb{N}+1-\beta$.  Unfortunately, there seem to
be no parameter choices for which the orthogonal polynomials on the
bi-lattice turn out to be in a well-known family and for which the
recurrence coefficients $a_n,b_n$ are explicitly known (as was the
case for $\beta=1/2, c_1=c_2$ in the generalized Charlier case).
Nevertheless, some parameter choices deserve special attention.  The
case $\beta=\gamma$ gives rise to the weight function
$w(x)=a^x/\Gamma(x+1)$.  With the choice $t=0$ this gives us the
Charlier polynomials on the lattice $\mathbb{N}$.  The case
$\beta=1/2, t=1$ gives the equally spaced lattice
$\frac{1}{2}\mathbb{N}$.  The weight in the lattice point $k/2$ is
then given by
\[w(k/2)=\frac{\Gamma(\gamma+k/2)}{\Gamma(\gamma)}\frac{(2\sqrt{a})^k}{k!},\]
as can be seen by using Legendre's duplication formula \cite[Eq.~5.5.5]{NIST}.

\section{Concluding remarks}
\begin{figure}[ht]
\centering
\resizebox{3in}{!}{\includegraphics{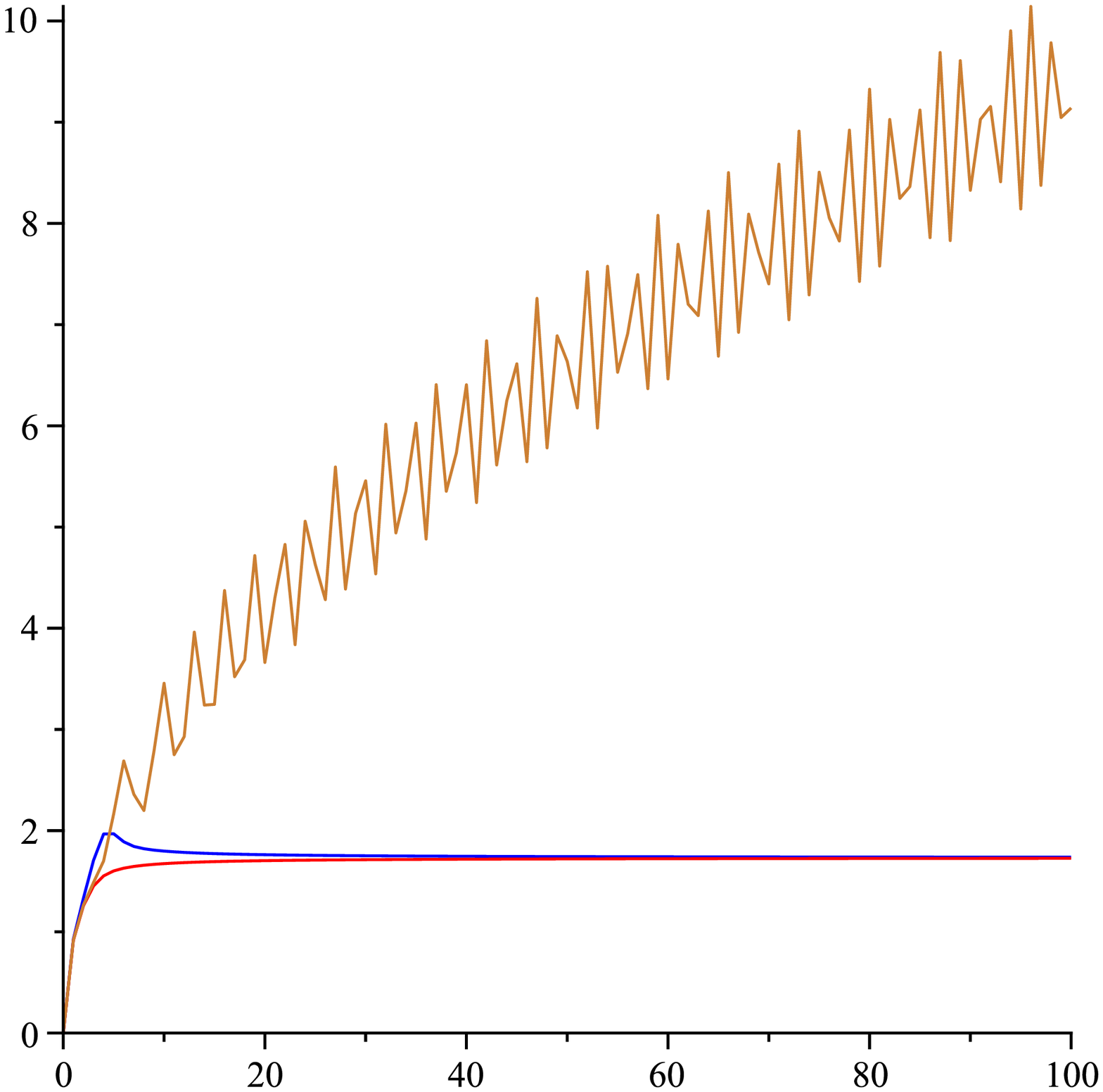}}
\resizebox{3in}{!}{\includegraphics{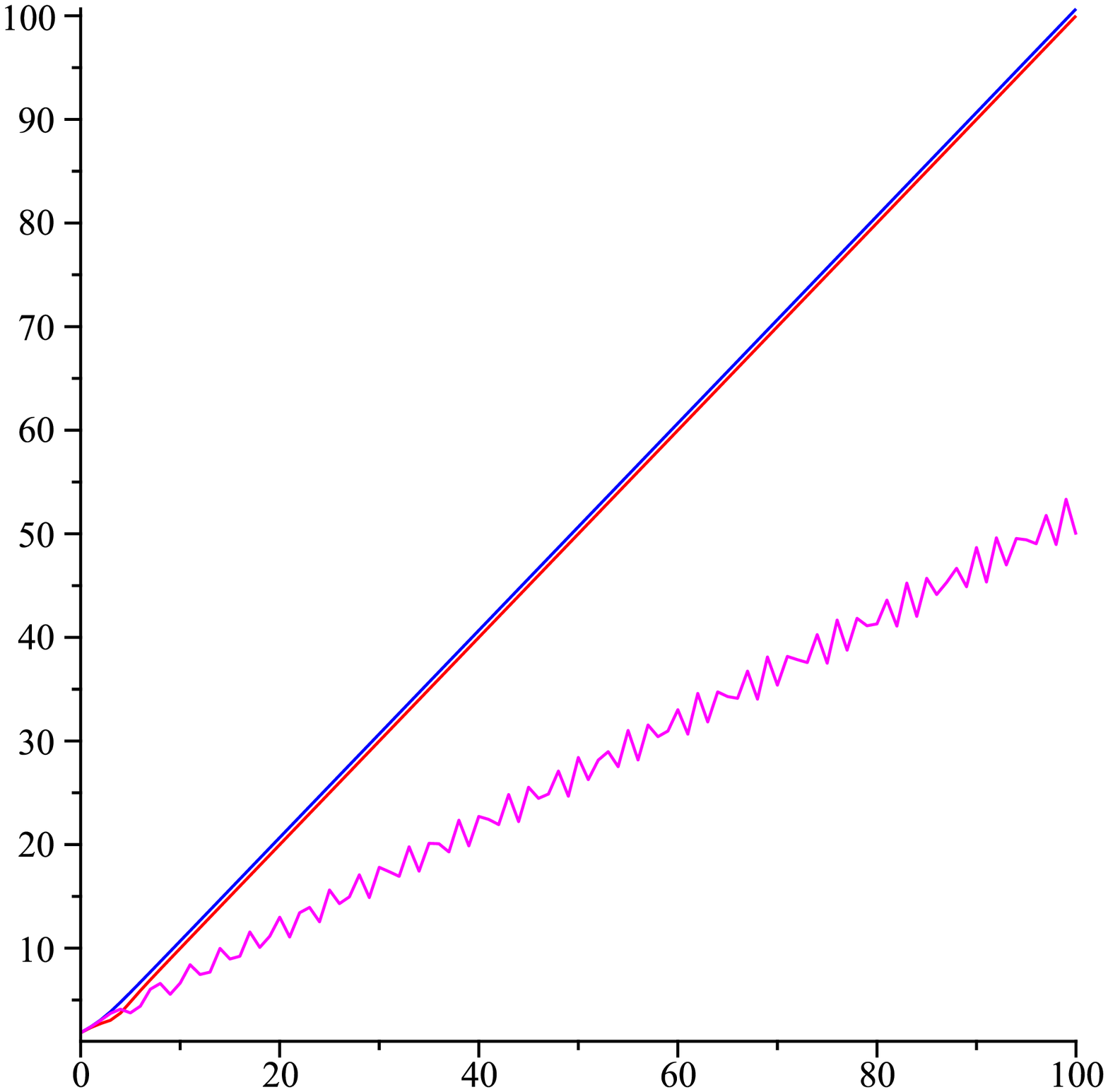}}
\caption{Recurrence coefficients for generalized Charlier polynomials ($a=3$, $\beta=1/3$, $t=10$)}
\label{fig:1}
\end{figure}

Even though the recurrence coefficients of the generalized Charlier
polynomials satisfy the same non-linear system of recurrence
equations, their behavior for the bi-lattice is quite different from
the behavior on the lattices. In Figure \ref{fig:1} we have plotted
on the left the recurrence coefficients $a_n$ and on the right the
$b_n$ for three cases. For the $a_n$ (left plot) the lowest curve
corresponds to the lattice $\mathbb{N}+1-\beta$, the curve just
above to the lattice $\mathbb{N}$ and the top curve corresponds to
the bi-lattice with $t=10$ (the parameters are $a=3$, $\beta=1/3$).
For the $b_n$ (right plot) the top curve corresponds to the lattice
$\mathbb{N}+1-\beta$, the curve just below to the lattice
$\mathbb{N}$ and the lowest curve to the bi-lattice with $t=10$.

Observe that the behavior for each of the two lattices is quite
similar and monotonic as $n \to \infty$, but for the bi-lattice the
behavior is oscillating with a different asymptotic behavior as $n
\to \infty$. We conjecture that for the lattice $\mathbb{N}$ and $\mathbb{N}+1-\beta$ one has the
asymptotic behavior $\lim_{n \to \infty} a_n^2 = a$, but for the bi-lattice one has
$a_n^2 = n\sqrt{a}/2 + \mathcal{O}(1)$, where the term $\mathcal{O}(1)$ is bounded but oscillating.
This behavior is in agreement with the special solution (\ref{eq:beta1/2}) for $\beta=1/2$. For the
$b_n$ we conjecture that
\[   \lim_{n \to \infty} b_n - n = \begin{cases} 
					0 & \textrm{for the lattice $\mathbb{N}$}, \\
					1-\beta & \textrm{for the lattice $\mathbb{N}+1-\beta$},
                                   \end{cases}  \]
and $b_n = n/2 + \mathcal{O}(1)$ for the bi-lattice. 

\begin{figure}[ht]
\centering
\resizebox{3in}{!}{\includegraphics{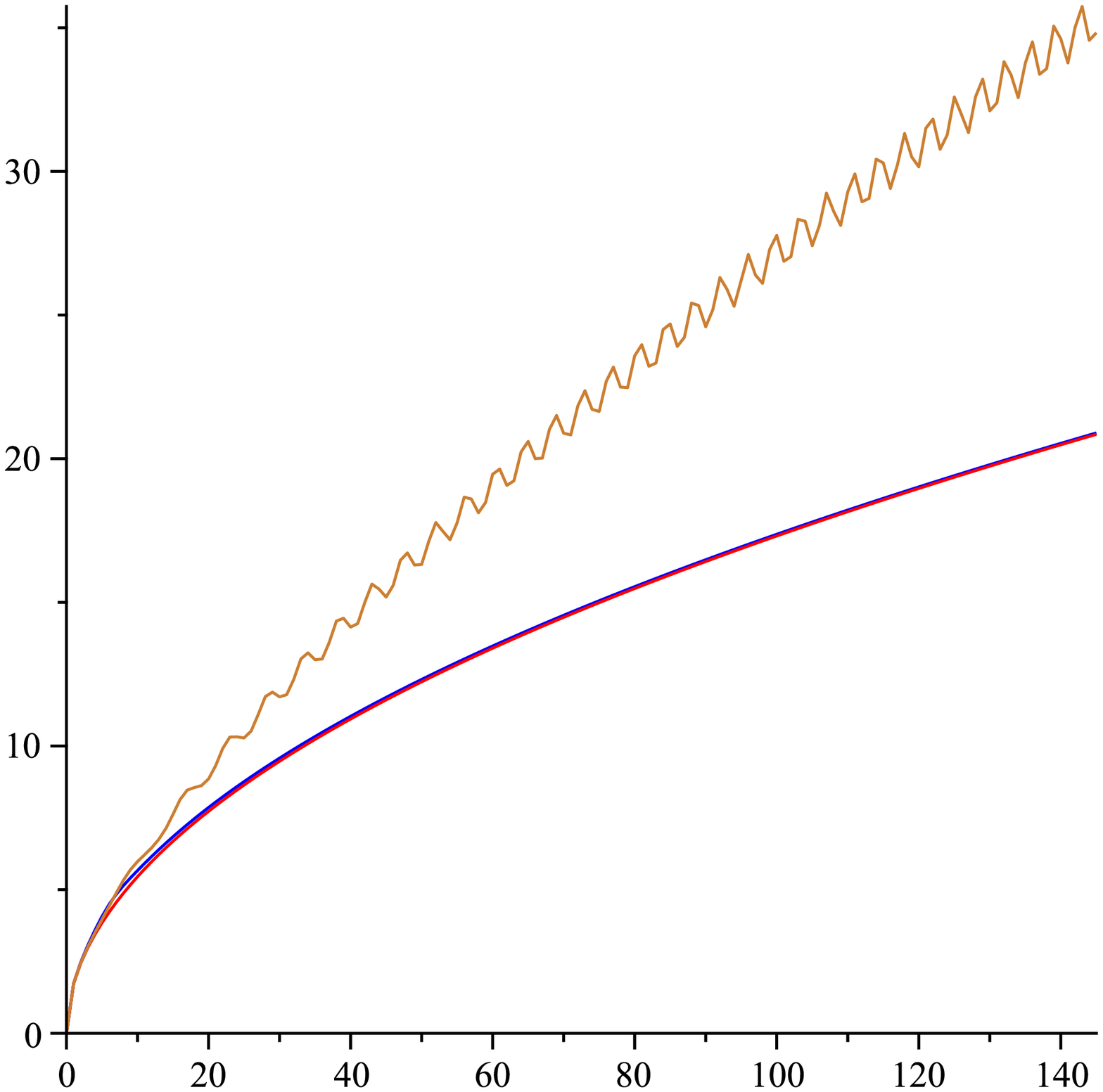}}
\resizebox{3in}{!}{\includegraphics{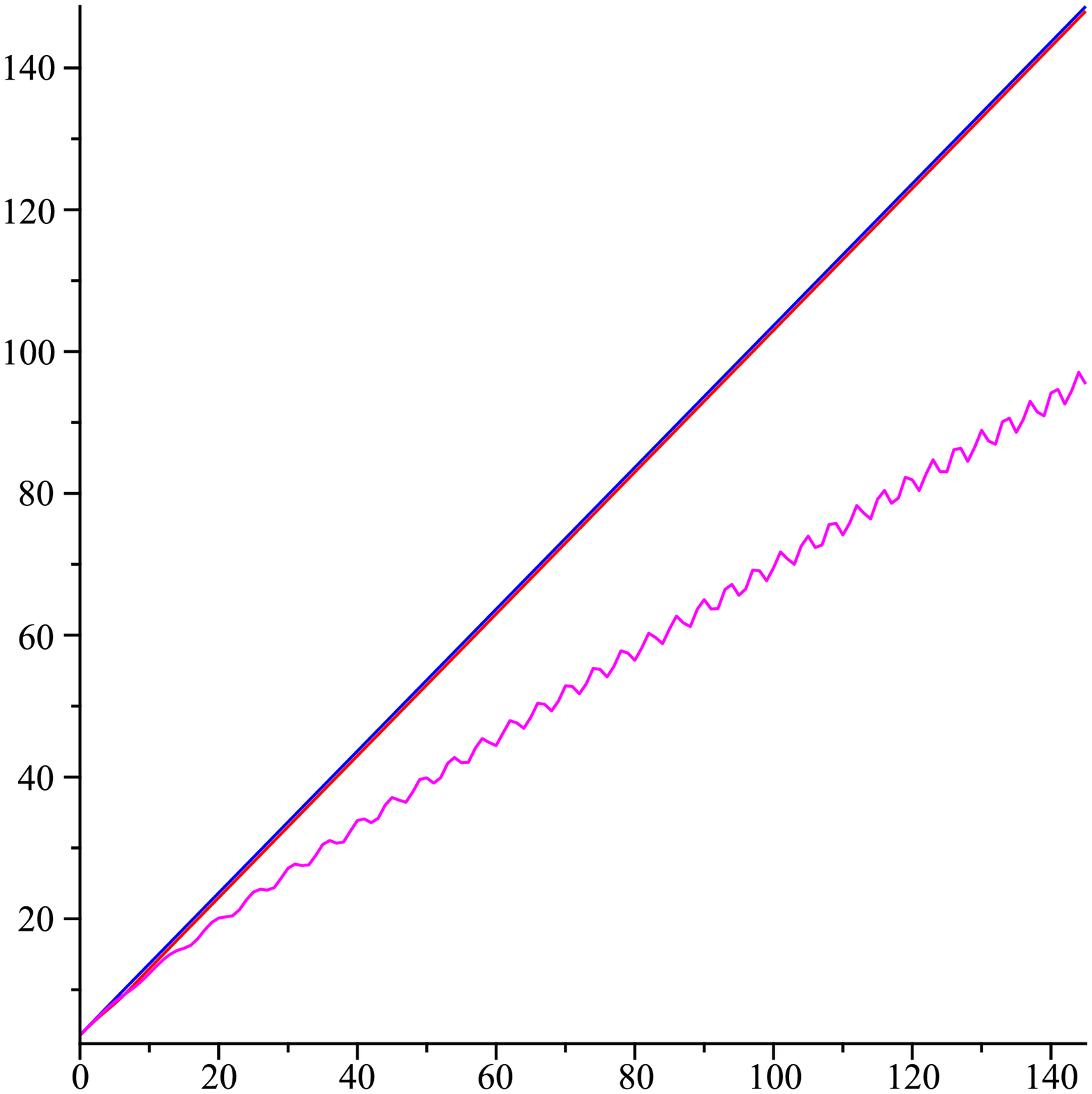}}
\caption{Recurrence coefficients for generalized Meixner polynomials ($a=3$, $\beta=2/3$, $\gamma = 9/10$, $t=2$)}
\label{fig:2}
\end{figure}
Similar observations hold for the recurrence coefficients of the
generalized Meixner polynomials. In Figure~\ref{fig:2} we have
plotted the recurrence coefficients $a_n$ (on the left) and $b_n$
(on the right) for the same three cases for the parameter values
$a=3$, $\beta=2/3$, $\gamma=9/10$. For the $a_n$, the lowest curve
corresponds to the lattice $\mathbb{N}+1-\beta$, the curve just
above to the lattice $\mathbb{N}$ (both curves are almost identical), and the top curve corresponds to
the bi-lattice with $t=2$.  For the $b_n$ the top curve corresponds
to the lattice $\mathbb{N}+1-\beta$, the one just below to the
lattice $\mathbb{N}$ (both curves are almost identical), and the lowest curve corresponds to the
bi-lattice. Here we conjecture that for each lattice one has
\[   \lim_{n \to \infty} a_n^2 - an = \begin{cases}
    (\gamma-\beta)a & \textrm{for the lattice $\mathbb{N}$}, \\
    (\gamma-1)a & \textrm{for the lattice $\mathbb{N}+1-\beta$,}
     \end{cases} \]
and
\[     \lim_{n \to \infty} b_n - n = \begin{cases}
     a & \textrm{for the lattice $\mathbb{N}$}, \\
     a + 1-\beta & \textrm{for the lattice $\mathbb{N}+1-\beta$.}
     \end{cases}  \]
The asymptotic behavior for
the bi-lattice is more difficult and we conjecture $a_n^2/n^{3/2} = \mathcal{O}(1)$ and $b_n/n = \mathcal{O}(1)$, where
the $\mathcal{O}(1)$ terms are oscillatory.

\begin{verbatim}
Christophe Smet (christophe@wis.kuleuven.be)
Walter Van Assche (walter@wis.kuleuven.be)
Department of Mathematics
Katholieke Universiteit Leuven
Celestijnenlaan 200B box 2400
BE-3001 Leuven, BELGIUM
\end{verbatim}

\end{document}